\documentclass[12pt]{article}

\usepackage{amsmath,amssymb}
\usepackage{tikz}
\usepackage{color}
\usepackage[toc]{appendix}
\usepackage{graphicx}
\usepackage{fancyhdr}
\usepackage{enumitem}
\usepackage{bbm}
\usepackage{graphicx}
\usepackage{parskip}
\usepackage{float}
\usepackage{chngpage}
\usepackage{calc}
\usepackage{bigints}
\usepackage{tikz}
\usepackage{array}
\usepackage{bbm}
\usepackage{booktabs}
\usepackage{rotating}
\usepackage{multirow}
\usepackage{adjustbox}
\usepackage{tabularx}
\usepackage{enumitem}
\usepackage{verbatim}
\usepackage{mathtools}
\usepackage{ragged2e}
\usepackage[makeroom]{cancel}
\usepackage{enumitem}
\usepackage{caption}
\usepackage{mathtools}
\usepackage{hyperref}
\usepackage{bbm}
\usepackage{amsmath,amsfonts,amsthm,bm}
\usepackage{amsmath,amssymb}
\usepackage{appendix}
\usepackage{graphicx}
\usepackage{pgfplots}
\usepackage{tikz}
\usepackage{verbatim}
\usepackage{amsthm}
\usepackage[english]{babel}
\usepackage{hyphenat}
\usepackage[makeindex]{imakeidx}
\usepackage{tikz}

\usetikzlibrary{datavisualization}
\usetikzlibrary{matrix}
\usetikzlibrary{datavisualization.formats.functions}

\newtheorem{theorem}{Theorem}[section]
\newtheorem{proposition}[theorem]{Proposition}

\newtheorem{remark}[theorem]{Remark}

\makeatletter
\def\section{\@startsection {section}{1}{\z@}{3.25ex plus 1ex minus
		.2ex}{1.5ex plus .2ex}{\large\bf}}
\def\subsection{\@startsection{subsection}{2}{\z@}{3.25ex plus 1ex minus
		.2ex}{1.5ex plus .2ex}{\normalsize\bf}}
\@addtoreset{equation}{section} 
\makeatother

\title{On a class of stochastic hyperbolic equations with double characteristics}

\author{Enrico Bernardi\thanks{Dipartimento di Scienze Statistiche Paolo Fortunati, Università di Bologna, Bologna, Italy. \textbf{e-mail}: enrico.bernardi@unibo.it} \and Alberto Lanconelli\thanks{Dipartimento di Scienze Statistiche Paolo Fortunati, Università di Bologna, Bologna, Italy. \textbf{e-mail}: alberto.lanconelli2@unibo.it}}

\date{\today}

\begin{document}
	
	\maketitle
	
	\bigskip
	
	\begin{abstract}
		We study the effect of Gaussian perturbations on a hyperbolic partial differential equation with double characteristics in two spatial dimensions. The coefficients of our partial differential operator depend polynomially on the space variables, while the noise is additive, white in time and coloured in space. We provide a sufficient condition on the spectral measure of the covariance functional describing the noise that allows for the existence of a random field solution for the resulting stochastic partial differential equation. Our approach is based on explicit computations for the fundamental solution of the partial differential operator and its Fourier transform.    
	\end{abstract}
	
	
	
	\bigskip
	
\allowdisplaybreaks	
	
\section{Introduction and statement of the main result}\label{intro}

The aim of this note is to investigate the stochastic linear hyperbolic equation
\begin{align}\label{SPDE usual}
\begin{cases}
\left(\partial_t^2-2\partial_t\partial_{x_1}-x_1^2\partial_{x_2}^2\right)u(t,x_1,x_2)=\dot{F}(t,x_1,x_2),&\quad t>0, (x_1,x_2)\in\mathbb{R}^2;\\
u(0,x_1,x_2)=0,&\quad(x_1,x_2)\in\mathbb{R}^2;\\
\partial_tu(0,x_1,x_2)=0,&\quad(x_1,x_2)\in\mathbb{R}^2,
\end{cases}
\end{align}
where 
\begin{align*}
F(\varphi):=\int_{\mathbb{R}^3}\varphi(t,x_1,x_2)\dot{F}(t,x_1,x_2)dtdx_1dx_2,\quad \varphi\in C_0^{\infty}(\mathbb{R}^3)
\end{align*}
is a family of Gaussian random variables, defined on a common complete probability space $(\Omega,\mathcal{A},\mathbb{P})$, with mean zero and covariance
\begin{align}\label{covariance preliminaries}
\mathbb{E}[F(\varphi)F(\psi)]=\int_0^{+\infty}\int_{\mathbb{R}^2}dx\int_{\mathbb{R}^2}dy\quad \!\!\!\!\varphi(t,x)f(x-y)\psi(t,y),\quad \varphi,\psi\in C_0^{\infty}(\mathbb{R}^3). 
\end{align}
The function $f:\mathbb{R}^2\to\mathbb{R}$ is assumed to be continuous in $\mathbb{R}^2\setminus\{0\}$ with $f(-x)=f(x)$, for all $x\in\mathbb{R}^2$.

The most relevant feature of the differential operator appearing in \eqref{SPDE usual} is the fact that its principal symbol, i.e.
\begin{align*}
p(x,\xi)=-\xi_0^2-2\xi_0\xi_1+x_1^2\xi_2^2,
\end{align*}
is hyperbolic with respect to $\xi_0$ and has double characteristics on the manifold
\begin{align*}
\Sigma=\{(x,\xi)\in \dot{T}^{\star}\mathbb{R}^3: x_1=0,\xi_0=0,\xi_1=0\},
\end{align*}
where $\dot{T}^{\star}\mathbb{R}^3$ denotes the phase-space cotangent bundle of $\mathbb{R}^3$ minus the $0$-section (see e.g. \cite{Hormander 1977},\cite{Hormander}). The fundamental matrix associated to
$p$ at a double point $\rho\in\Sigma$ is then computed as $F_p(\rho)=\frac{1}{2}dH_p(\rho)$,
where $H_p$ is the Hamilton vector field of $p$ and it is readily seen that $F_p$ has a Jordan block of order $4$ at the eigenvalue $0$; this is one of the (well-known) three possible non-effectively hyperbolic cases, if one excludes the case of a pair of non zero-real eigenvalues (the so called effectively hyperbolic case). \\
As it is the case with hyperbolic operators with multiple characteristics, lower order terms may modify the behaviour of the well-posedness of the Cauchy problem. A large number of papers has been devoted to that question in the deterministic setting and in our case when studying the problem
\begin{align*}
(\partial_t^2-2\partial_t\partial_{x_1}-x_1^2\partial_{x_2}^2+\mathtt{s}\partial_{x_2})u(t,x_1,x_2)=g(t,x_1,x_2)
\end{align*}
the behaviour changes whether $\mathtt{s}=0$ or not. This is a special case of the Ivrii-Petkov conditions (see \cite{IP}) and we will deal with that and its effects on random perturbations in a following paper.

We recall that the analogue of (\ref{SPDE usual}) in the case of strictly hyperbolic operators, or the wave equation for that matter, has been thoroughly studied in the literature for different spatial dimensions. We mention for instance \cite{Orsingher}, \cite{DF}, \cite{DM}, \cite{Peszat}, \cite{Sanz-Sole} and the reference quoted there. The framework adopted in this paper is the one proposed in \cite{Dalang 99} which extends the classical set up of \cite{Walsh}. We also mention \cite{DQ} for a comparison between the abstract framework of \cite{Da Prato} and the one discussed in \cite{Dalang 99}. Lastly, we mention that the analysis of some non strictly hyperbolic stochastic partial differential equations has been carried in the recent papers \cite{Ascanelli 2019} and \cite{Ascanelli 2018}. There the authors adopt the stochastic framework proposed in \cite{Dalang 99} and prove existence of random field solutions; however, their class of operators does not cover the one treated here.

To state our main theorem, we now shortly describe the framework, referring the reader to \cite{Walsh} and \cite{Dalang 99} for further details.\\
We denote by $\mathcal{D}(\mathbb{R}^3)$ the space of functions $\varphi\in C_0^{\infty}(\mathbb{R}^3)$ endowed with the topology induced by the following notion of convergence: $\varphi_n\to\varphi$ if 
\begin{enumerate}
	\item there exists a compact set $K$ of $\mathbb{R}^3$ such that the support of $\varphi_n-\varphi$ is contained in $K$, for all $n\geq 1$;
	\item $\lim_{n\to+\infty}D^{\alpha}\varphi_n=D^{\alpha}\varphi$, uniformly on $K$ for each multiindex $\alpha$.
\end{enumerate}
A direct verification using identity (\ref{covariance preliminaries}) shows that the map $\varphi\mapsto F(\varphi)$ is linear and continuous in $\mathbb{L}^2(\Omega)$; this implies that $F$ has a version with values in $\mathcal{D}'(\mathbb{R}^3)$ which in turn allows for a distributional $\omega$-wise interpretation of the partial differential equation \eqref{SPDE usual}. For the distributional solution to be a real valued measurable stochastic process, we need to extend $F$ to a worthy martingale and interpret the distributional solution 
\begin{align*}
u(t,x_1,x_2)=\int_0^t\int_{\mathbb{R}^2}\Gamma(t-s,x_1,y_1,x_2-y_2)\dot{F}(s,y_1,y_2)dy_1dy_2,
\end{align*}
as a stochastic integral; here $\Gamma$ denotes the fundamental solution of the differential operator $\partial_t^2-2\partial_t\partial_{x_1}-x_1^2\partial_{x_2}^2$ from \eqref{SPDE usual} (see Section \ref{fundamental solution section} below). To this aim, by suitably approximating indicator functions of bounded Borel subsets of $[0,+\infty[\times\mathbb{R}^2$ with elements from $\mathcal{D}(\mathbb{R}^3)$ and employing the $\mathbb{L}^2(\Omega)$-continuity mentioned above, we first define $F(A):=F\left(\boldsymbol{1}_A\right)$, $A\in\mathcal{B}_b([0,+\infty[\times\mathbb{R}^2)$ and 
\begin{align*}
M_t(B):=F([0,t]\times B),\quad t\geq 0, B\in\mathcal{B}_b(\mathbb{R}^2).
\end{align*}
Then, if we let
\begin{align*}
\mathcal{F}_t^0:=\sigma(M_s(B), 0\leq s\leq t, B\in\mathcal{B}_b(\mathbb{R}^2)),\quad \mathcal{F}_t:=\mathcal{F}_t^0\vee\mathcal{N},
\end{align*}
where $\mathcal{N}$ denotes the $\sigma$-algebra generated by $\mathbb{P}$-null sets, we get that
\begin{align*}
\left(\{M_t(B)\}_{t\geq 0,B\in\mathcal{B}_b(\mathbb{R}^2)},\{\mathcal{F}_t\}_{t\geq 0}\right)
\end{align*}
is a \emph{worthy martingale measure}. By construction, for all $B\in\mathcal{B}_b(\mathbb{R}^2)$ the stochastic process $\{M_t(B)\}_{t\geq 0}$ is a continuous martingale and we have
\begin{align}\label{worthy def}
F(\varphi)=\int_0^{+\infty}\int_{\mathbb{R}^2}\varphi(t,x)M(dt,dx).
\end{align}  
On the other hand, using elementary properties of the Fourier transform, we can rewrite identity (\ref{covariance preliminaries}) as
\begin{align}\label{FT covariance}
\mathbb{E}[F(\varphi)F(\psi)]=\int_0^{+\infty}\int_{\mathbb{R}^2}\mathcal{F}\varphi(t,\xi)\overline{\mathcal{F}\psi(t,\xi)}d\mu(\xi)dt,
\end{align}
where $\mathcal{F}\eta$ denotes the Fourier transform of $\eta$, i.e.
\begin{align*}
\mathcal{F}\eta(\xi):=\int_{\mathbb{R}^2}e^{-i\xi\cdot x}\eta(x)dx,\quad \xi\in\mathbb{R}^2,
\end{align*}
and $\mu$, the \emph{spectral measure of $f$}, is a non-negative tempered measure $\mu$ on $\mathbb{R}^2$ such that
\begin{align*}
\int_{\mathbb{R}^2}f(x)\eta(x)dx=\int_{\mathbb{R}^2}\mathcal{F}\eta(\xi)d\mu(\xi),\quad\mbox{ for all $\eta\in S(\mathbb{R}^2)$}.
\end{align*} 
Combining identity \eqref{FT covariance} with \eqref{worthy def} we get
\begin{align}\label{isometry}
\mathbb{E}\left[\left|\int_0^{+\infty}\int_{\mathbb{R}^2}\varphi(t,x)M(dt,dx)\right|^2\right]=\int_0^{+\infty}\int_{\mathbb{R}^2}|\mathcal{F}\varphi(t,\xi)|^2d\mu(\xi)dt.
\end{align}
The last isometry determines the class of admissible deterministic integrands for the stochastic integral in (\ref{worthy def}). We will say that $\{u(t,x_1,x_2)\}_{t\geq 0, (x_1,x_2)\in\mathbb{R}^2}$ is a \emph{random field solution} to (\ref{SPDE usual}) if
 \begin{align}\label{worthy intro}
u(t,x_1,x_2):=\int_0^t\int_{\mathbb{R}^2}\Gamma(t-s,x_1,y_1,x_2-y_2)M(ds,dy_1,dy_2),
\end{align}
is a well defined stochastic integral (i.e. the right hand side in \eqref{isometry} is finite) and 
the map
\begin{align*}
[0,+\infty[\times\mathbb{R}^2\ni(t,x_1,x_2)\mapsto u(t,x_1,x_2)
\end{align*} 
is measurable. We are now ready to state the main theorem of the present paper; the proof can be found in Sections 2 and 3 (closed-form expression for the fundamental solution and existence for the random field solution, respectively).

\begin{theorem}\label{main theorem}
Assume the spectral measure $\mu$ to satisfy
\begin{align}\label{SC}
\int_{\mathbb{R}^2}\frac{1}{1+|\xi|^{2/3}}d\mu(\xi)<+\infty.
\end{align}
Then, the stochastic partial differential equation \eqref{SPDE usual} admits
a random field solution $\{u(t,x)\}_{t\in [0,T], x\in\mathbb{R}^2}$ with representation \eqref{worthy intro} where
\begin{align*}
\Gamma(t,x_1,y_1,x_2)=
\begin{cases}
\frac{\sqrt{3}}{2\pi}\frac{1}{\sqrt{(y_1^3-x_1^3)(2t+x_1-y_1)-3x_2^2}},&\mbox{if }(t,x_1,x_2)\in A_{t,x_1,x_2},\\
0,&\mbox{ if }(t,x_1,x_2)\notin A_{t,x_1,x_2},
\end{cases}
\end{align*}
and 
\begin{align*}
A_{t,x_1,x_2}:=\left\{(t,x_1,x_2)\in\mathbb{R}^3: t>0,y_1>x_1,(y_1^3-x_1^3)(2t+x_1-y_1)-3x_2^2>0\right\}.
\end{align*}
\end{theorem}

\begin{remark}
It is proved in \cite{Dalang 99} that for the stochastic wave equation
	\begin{align*}
	\begin{cases}
	\left(\partial_t^2-\Delta\right)u(t,x)=\dot{F}(t,x),&\quad t>0, x\in\mathbb{R}^d;\\
	u(0,x)=0,&\quad x\in\mathbb{R}^d;\\
	\partial_tu(0,x)=0,&\quad x\in\mathbb{R}^d,
	\end{cases}
	\end{align*}
with $\Delta:=\partial_{x_1}^2+\cdot\cdot\cdot+\partial_{x_d}^2$, a sufficient condition for the existence of random field solutions is
\begin{align*}
\int_{\mathbb{R}^d}\frac{1}{1+|\xi|^{2}}d\mu(\xi)<+\infty,
\end{align*}
for any spatial dimension $d\geq 1$. A comparison with (\ref{SC}) shows that the existence of a  random filed solution for (\ref{SPDE usual}) requires more stringent assumptions on the spectral measure $\mu$, and hence on the Gaussian noise $F$, than its simplest strictly hyperbolic counterpart.
\end{remark}

The paper is organized as follows: in Section 2 we derive an explicit representation for the fundamental solution of the partial differential operator in \eqref{SPDE usual} while in Section 3 we prove the existence for a random field solution under the integrability condition \eqref{SC}.

\section{The fundamental solution}\label{fundamental solution section}

In this section we describe a derivation of the fundamental solution for the partial differential operator in (\ref{SPDE usual}). Following \cite{Hormander}, we set $D_t:=\frac{1}{i}\partial_{t}$, $D_{x_j}:=\frac{1}{i}\partial_{x_j}$, for $j=1,2$, and consider the problem
\begin{align}\label{PDE with delta}
\left(-D_t^2+2D_tD_{x_1}+x_1^2D_{x_2}^2\right)\Gamma(t,x_1,y_1,x_2)=\delta_{(0,y_1,0)}(t,x_1,x_2),
\end{align}
for $t>0$ and $(x_1,x_2)\in\mathbb{R}^2$; here $\delta_z$ stands for the Dirac's delta distribution with mass at $z\in\mathbb{R}^3$ and $y_1\in\mathbb{R}$ is a fixed parameter. We observe that the operator under investigation is not invariant by translation in the variable $x_1$; hence, the parameter $y_1$ serves to keep trace of this fact. \\
We now denote by $\hat{\Gamma}(t,x_1,y_1,\xi_2)$ the Fourier transform w.r.t. $x_2$ of $x_2\mapsto \Gamma(t,x_1,y_1,x_2)$, i.e.
\begin{align*}
\hat{\Gamma}(t,x_1,y_1,\xi_2):=\int_{\mathbb{R}}e^{-i\xi_2x_2}\Gamma(t,x_1,y_1,x_2)dx_2,\quad \xi_2\in\mathbb{R},
\end{align*}
and by $\Gamma^{\dagger}(\xi_0,x_1,y_1,x_2)$ the Fourier-Laplace transform w.r.t. $t$ of $t\mapsto \Gamma(t,x_1,y_1,x_2)$, i.e.
\begin{align*}
\Gamma^{\dagger}(\xi_0,x_1,y_1,x_2):=\int_{0}^{+\infty}e^{-i\xi_0t}\Gamma(t,x_1,y_1,x_2)dt,\quad \mathtt{Im}(\xi_0)<0.
\end{align*} 
Transforming equation (\ref{PDE with delta}) we get
\begin{align}\label{a}
\left(-\xi_0^2+2\xi_0D_{x_1}+x_1^2\xi_2^2\right)\hat{\Gamma}^{\dagger}(\xi_0,x_1,y_1,\xi_2)=\delta(x_1-y_1),
\end{align}
with $\delta$ being now the Dirac's delta distribution with mass at $0\in\mathbb{R}$. Recalling that $D_{x_1}=\frac{1}{i}\partial_{x_1}$ and setting
\begin{align*}
q(x_1;\xi_0,\xi_2):=x_1^2\xi_2^2-\xi_0^2,
\end{align*}
we can rewrite equation (\ref{a}) as
\begin{align}\label{b}
\left(\partial_{x_1}+\frac{i}{2\xi_0}q(x_1;\xi_0,\xi_2)\right)\hat{\Gamma}^{\dagger}(\xi_0,x_1,y_1,\xi_2)=\frac{i}{2\xi_0}\delta(x_1-y_1).
\end{align}
The homogeneous part of the last equation can be solved as
\begin{align*}
v(\xi_0,x_1,\xi_2)&=C(\xi_0,\xi_2)e^{-\frac{i}{2\xi_0}\int_0^{x_1}q(t;\xi_0,\xi_2)dt}\\
&=C(\xi_0,\xi_2)e^{-\frac{i}{2\xi_0}\int_0^{x_1}t^2\xi_2^2-\xi_0^2dt}\\
&=C(\xi_0,\xi_2)e^{-\frac{i}{2\xi_0}\left(\frac{\xi_2^2}{3}x_1^3-\xi_0^2x_1\right)}\\
&=C(\xi_0,\xi_2)e^{i\frac{\xi_0}{2}x_1-i\frac{\xi_2^2}{6\xi_0}x_1^3}.
\end{align*}
On the other hand, since it is natural to expect that $\hat{\Gamma}^{\dagger}(\xi_0,x_1,y_1,\xi_2)$ vanishes when $x_1>y_1$, we choose
\begin{align*}
\hat{\Gamma}^{\dagger}(\xi_0,x_1,y_1,\xi_2)=v(\xi_0,x_1,\xi_2)H(y_1-x_1),
\end{align*} 
as a particular solution of (\ref{b}) (here $x\mapsto H(x)$ denotes the Heaviside function). We now have to find the constant $C(\xi_0,\xi_2)$ that makes identity (\ref{b}) true:
\begin{align*}
&\left(\partial_{x_1}+\frac{i}{2\xi_0}q(x_1;\xi_0,\xi_2)\right)\hat{\Gamma}^{\dagger}(\xi_0,x_1,y_1,\xi_2)\\
&\quad=\partial_{x_1}\hat{\Gamma}^{\dagger}(\xi_0,x_1,y_1,\xi_2)+\frac{i}{2\xi_0}q(x_1;\xi_0,\xi_2)\hat{\Gamma}^{\dagger}(\xi_0,x_1,y_1,\xi_2)\\
&\quad=\partial_{x_1}(v(\xi_0,x_1,\xi_2)H(y_1-x_1))+\frac{i}{2\xi_0}q(x_1;\xi_0,\xi_2)\hat{\Gamma}^{\dagger}(\xi_0,x_1,y_1,\xi_2)\\
&\quad=(\partial_{x_1}v(\xi_0,x_1,\xi_2))H(y_1-x_1)+v(\xi_0,x_1,\xi_2)\partial_{x_1}H(y_1-x_1)\\
&\quad\quad+\frac{i}{2\xi_0}q(x_1;\xi_0,\xi_2)\hat{\Gamma}^{\dagger}(\xi_0,x_1,y_1,\xi_2)\\
&\quad=-\frac{i}{2\xi_0}q(x_1;\xi_0,\xi_2)v(\xi_0,x_1,\xi_2)H(y_1-x_1)-v(\xi_0,x_1,\xi_2)\delta(y_1-x_1)\\
&\quad\quad+\frac{i}{2\xi_0}q(x_1;\xi_0,\xi_2)\hat{\Gamma}^{\dagger}(\xi_0,x_1,y_1,\xi_2)\\
&\quad=-v(\xi_0,x_1,\xi_2)\delta(y_1-x_1)\\
&\quad=-C(\xi_0,\xi_2)e^{i\frac{\xi_0}{2}x_1-i\frac{\xi_2^2}{6\xi_0}x_1^3}\delta(y_1-x_1).
\end{align*}
This gives
\begin{align*}
C(\xi_0,\xi_2)=-\frac{i}{2\xi_0}e^{-i\frac{\xi_0}{2}y_1+i\frac{\xi_2^2}{6\xi_0}y_1^3},
\end{align*}
and hence
\begin{align*}
\hat{\Gamma}^{\dagger}(\xi_0,x_1,y_1,\xi_2)=-\frac{i}{2\xi_0}e^{i\frac{\xi_0}{2}(x_1-y_1)-i\frac{\xi_2^2}{6\xi_0}(x_1^3-y_1^3)}H(y_1-x_1).
\end{align*} 
We now proceed inverting the transforms; to this aim we write the last expression as
\begin{align*}
\hat{\Gamma}^{\dagger}(\xi_0,x_1,y_1,\xi_2)=G(\xi_0,x_1)e^{-\frac{A(\xi_0)}{2}\xi_2^2}, 
\end{align*} 
where
\begin{align}\label{c1}
G(\xi_0,x_1,y_1):=-\frac{i}{2\xi_0}e^{i\frac{\xi_0}{2}(x_1-y_1)}H(y_1-x_1),
\end{align}    
and
\begin{align}\label{c2}
A(\xi_0,x_1,y_1):=i\frac{x_1^3-y_1^3}{3\xi_0}.
\end{align}
Then,
\begin{align}\label{c}
\Gamma^{\dagger}(\xi_0,x_1,y_1,x_2)&=\frac{1}{2\pi}\int_{\mathbb{R}}e^{ix_2\xi_2}\hat{\Gamma}^{\dagger}(\xi_0,x_1,y_1,\xi_2)d\xi_2\nonumber\\
&=\frac{1}{2\pi}\int_{\mathbb{R}}e^{ix_2\xi_2}G(\xi_0,x_1,y_1)e^{-\frac{A(\xi_0,x_1,y_1)}{2}\xi_2^2}d\xi_2\nonumber\\
&=\frac{1}{2\pi}G(\xi_0,x_1,y_1)\int_{\mathbb{R}}e^{ix_2\xi_2-\frac{A(\xi_0,x_1,y_1)}{2}\xi_2^2}d\xi_2\nonumber\\
&=\frac{1}{\sqrt{2\pi}}G(\xi_0,x_1,y_1)A(\xi_0,x_1,y_1)^{-\frac{1}{2}}e^{-\frac{x_2^2}{2}A(\xi_0,x_1,y_1)^{-1}}.
\end{align}
In the last equality we utilized the identity
\begin{align*}
\int_{\mathbb{R}}e^{-ixy-\frac{1}{2}Ax^2}dx=\sqrt{2\pi}A^{-\frac{1}{2}}e^{-\frac{1}{2}A^{-1}y^2},
\end{align*}
with $A^{\frac{1}{2}}$ chosen in a such a way that $A^{\frac{1}{2}}>0$, if $\mathtt{Re}(A)>0$. In our case, $A(\xi_0,x_1,y_1)=i\frac{x_1^3-y_1^3}{3\xi_0}$, $\mathtt{Im}(\xi_0)<0$ and $x_1<y_1$; this means that $\xi_0=|\xi_0|e^{i\theta}$, for some $\theta\in ]-\pi,0[$, that $x_1^3-y_1^3=|x_1^3-y_1^3|e^{-i\pi}$ and hence
\begin{align*}
A(\xi_0,x_1,y_1)=e^{i\frac{\pi}{2}}\frac{|x_1^3-y_1^3|e^{-i\pi}}{3|\xi_0|e^{i\theta}}=\frac{|x_1^3-y_1^3|}{3|\xi_0|}e^{i(-\frac{\pi}{2}-\theta)}.
\end{align*}
Since $-\frac{\pi}{2}-\theta\in ]-\frac{\pi}{2},\frac{\pi}{2}[$, we conclude that $\mathtt{Re}(A(\xi_0,x_1,y_1))>0$. Replacing in (\ref{c}) the definitions of $G(\xi_0,x_1,y_1)$ and $A(\xi_0,x_1,y_1)$ from (\ref{c1}) and (\ref{c2}), respectively, we get
\begin{align}\label{d}
\Gamma^{\dagger}(\xi_0,x_1,y_1,x_2)&=-\frac{1}{\sqrt{2\pi}}\frac{i}{2\xi_0}e^{i\frac{\xi_0}{2}(x_1-y_1)}
\left(\frac{y_1^3-x_1^3}{3i\xi_0}\right)^{-\frac{1}{2}}e^{\frac{i}{2}\frac{3\xi_0}{x_1^3-y_1^3}x_2^2}H(y_1-x_1)\nonumber\\
&=-\frac{1}{\sqrt{2\pi}}\frac{i}{2\xi_0}
\left(\frac{3i\xi_0}{y_1^3-x_1^3}\right)^{\frac{1}{2}}e^{i\frac{\xi_0}{2}\left(x_1-y_1+\frac{3}{x_1^3-y_1^3}x_2^2\right)}H(y_1-x_1)\nonumber\\
&=\frac{1}{\sqrt{2\pi}}\frac{1}{2}\frac{1}{i\xi_0}
\left(\frac{3}{y_1^3-x_1^3}\right)^{\frac{1}{2}}(i\xi_0)^{\frac{1}{2}}e^{i\frac{\xi_0}{2}\left(x_1-y_1+\frac{3}{x_1^3-y_1^3}x_2^2\right)}H(y_1-x_1)\nonumber\\
&=\frac{1}{\sqrt{2\pi}}\frac{1}{2}
\left(\frac{3}{y_1^3-x_1^3}\right)^{\frac{1}{2}}(i\xi_0)^{-\frac{1}{2}}e^{i\frac{\xi_0}{2}\left(x_1-y_1+\frac{3}{x_1^3-y_1^3}x_2^2\right)}H(y_1-x_1).
\end{align}
We now observe that the function $\xi_0\mapsto (i\xi_0)^{-\frac{1}{2}}$ is the Fourier-Laplace tranform of $t\mapsto (\pi t)^{-\frac{1}{2}}H(t)$. In fact, writing $\xi_0=\alpha_0+i\beta_0$ and recalling that $\beta_0<0$, we get
\begin{align*}
\int_0^{+\infty}e^{-i\xi_0t}(\pi t)^{-\frac{1}{2}}dt&=\int_0^{+\infty}e^{-i\alpha_0t+\beta_0t}(\pi t)^{-\frac{1}{2}}dt\\
&=\int_0^{+\infty}e^{-i\alpha_0t}\frac{1}{\sqrt{\pi}}t^{1/2-1}e^{-(-\beta_0)t}dt\\
&=(-\beta_0)^{-1/2}\int_0^{+\infty}e^{-i\alpha_0t}\frac{(-\beta_0)^{1/2}}{\sqrt{\pi}}t^{1/2-1}e^{-(-\beta_0)t}dt\\
&=(-\beta_0)^{-1/2}\left(1-i\frac{\alpha_0}{\beta_0}\right)^{-1/2}\\
&=(i\xi_0)^{-\frac{1}{2}}.
\end{align*}
In the fourth equality above we recognize the characteristic function of a Gamma distribution with parameters $(\frac{1}{2},-\beta_0)$. Lastly, denoting the function $t\mapsto (\pi t)^{-\frac{1}{2}}H(t)$ with $\chi(t)$ and $\left(x_1-y_1+\frac{3}{x_1^3-y_1^3}x_2^2\right)/2$ with $\mathcal{T}$, we can continue in \eqref{d} as
\begin{align*}
\Gamma^{\dagger}(\xi_0,x_1,y_1,x_2)&=\frac{1}{\sqrt{2\pi}}\frac{1}{2}
\left(\frac{3}{y_1^3-x_1^3}\right)^{\frac{1}{2}}\chi^{\dagger}(\xi_0)e^{i\xi_0\mathcal{T}}H(y_1-x_1)\\
&=\frac{1}{\sqrt{2\pi}}\frac{1}{2}
\left(\frac{3}{y_1^3-x_1^3}\right)^{\frac{1}{2}}(\chi(\cdot+\mathcal{T}))^{\dagger}(\xi_0)H(y_1-x_1).
\end{align*}
This entails
\begin{align*}
\Gamma(t,x_1,y_1,x_2)&=\frac{1}{\sqrt{2\pi}}\frac{1}{2}
\left(\frac{3}{y_1^3-x_1^3}\right)^{\frac{1}{2}}(\pi (t+\mathcal{T}))^{-\frac{1}{2}}H(t)H(t+\mathcal{T})H(y_1-x_1)\\
&=\frac{\sqrt{3}}{2\pi}\frac{1}{\sqrt{(y_1^3-x_1^3)(2t+x_1-y_1)-3x_2^2}}H(t)H(t+\mathcal{T})H(y_1-x_1).
\end{align*}
This completes the proof of the following result.

\begin{proposition}\label{FS prop}
The fundamental solution of the operator $\partial_t^2-2\partial_t\partial_{x_1}-x_1^2\partial_{x_2}^2$ on $]0,+\infty[\times\mathbb{R}^2$ is given by 
\begin{align*}
\Gamma(t,x_1,y_1,x_2)=
\begin{cases}
\frac{\sqrt{3}}{2\pi}\frac{1}{\sqrt{(y_1^3-x_1^3)(2t+x_1-y_1)-3x_2^2}},&\mbox{if }(t,x_1,x_2)\in A_{t,x_1,x_2},\\
0,&\mbox{ if }(t,x_1,x_2)\notin A_{t,x_1,x_2},
\end{cases}
\end{align*}
where
\begin{align*}
A_{t,x_1,x_2}:=\left\{(t,x_1,x_2)\in\mathbb{R}^3: t>0,y_1>x_1,(y_1^3-x_1^3)(2t+x_1-y_1)-3x_2^2>0\right\}.
\end{align*}
\end{proposition}

\section{Existence of the random field solution}\label{main section}

We now prove that
\begin{align}\label{worthy}
u(t,x_1,x_2):=\int_0^t\int_{\mathbb{R}^2}\Gamma(t-s,x_1,y_1,x_2-y_2)M(ds,dy_1,dy_2),
\end{align}
where $M$ denotes the worthy martingale measure associated with $F$ and $\Gamma$ the function from Proposition \ref{FS prop}, is a random field solution for the Cauchy problem \eqref{SPDE usual}. To do that we first need to verify the bound
\begin{align}\label{condition for integral}
\int_0^t\int_{\mathbb{R}^2}|\mathcal{F}\Gamma(t-s,x_1,\cdot,x_2-\cdot)(\xi_1,\xi_2)|^2\mu(d\xi_1d\xi_2)ds<+\infty,
\end{align}
for the Fourier transform of the function
\begin{align*}
(y_1,y_1)\mapsto \Gamma(t-s,x_1,y_1,x_2-y_2);
\end{align*}
this will ensure that $u(t,x_1,x_2)$ is a well defined stochastic integral. To ease the notation we set
\begin{align*}
h:=h(t-s,x_1,y_1):=(2(t-s)+x_1-y_1)(y_1^3-x_1^3)/3;
\end{align*}
then, 
\begin{align*}
&\mathcal{F}\Gamma(t-s,x_1,\cdot,x_2-\cdot)(\xi)\\
&\quad=\int_{\mathbb{R}}dy_1e^{-iy_1\xi_1}\int_{\mathbb{R}} dy_2e^{-iy_2\xi_2}\Gamma(t-s,x_1,\xi,x_2-\xi_2)\\
&\quad=\frac{1}{2\pi}\int dy_1e^{-iy_1\xi_1}\int dy_2e^{-iy_2\xi_2}\frac{1}{\sqrt{h-(x_2-y_2)^2}}H(h-(x_2-y_2)^2)H(y_1-x_1).
\end{align*}
Observe that condition $h(t-s,x_1,y_1)-(x_2-y_2)^2>0$ implies $h(t-s,x_1,y_1)>0$ which, in combination with $x_1<y_1$, gives $x_1<y_1<x_1+2(t-s)$. Therefore, making the change of variable $\sigma:=\frac{y_2-x_2}{\sqrt{h(t-s,x_1,y_1)}}$ and observing that
\begin{align*}
h(t-s,x_1,y_1)-(x_2-y_2)^2>0\quad\Longleftrightarrow\quad\sigma\in ]-1,1[,
\end{align*}
we get
\begin{align}\label{aa}
\mathcal{F}\Gamma(t-s,x_1,\cdot,x_2-\cdot)(\xi)&=\frac{1}{2\pi}\int_{x_1}^{x_1+2(t-s)} dy_1e^{-iy_1\xi_1}\int_{-1}^1 e^{-i(x_2+\sigma\sqrt{h})\xi_2}\frac{1}{\sqrt{1-\sigma^2}}d\sigma\nonumber\\
&=\frac{e^{-ix_2\xi_2}}{2\pi}\int_{x_1}^{x_1+2(t-s)} dy_1e^{-iy_1\xi_1}\int_{-\frac{\pi}{2}}^{\frac{\pi}{2}} e^{-i\sqrt{h}\sin(\theta)\xi_2}d\theta;
\end{align}
here we performed the further change of variable $\sigma:=\sin(\theta)$, $\theta\in ]-\pi/2,\pi/2[$. Note that the last integral above can be written as
\begin{align*}
\int_{-\frac{\pi}{2}}^{\frac{\pi}{2}} e^{-i\sqrt{h}\sin(\theta)\xi_2}d\theta&=\int_{-\frac{\pi}{2}}^{\frac{\pi}{2}} \cos\left(\sqrt{h}\sin(\theta)\xi_2\right)d\theta-i\int_{-\frac{\pi}{2}}^{\frac{\pi}{2}} \sin\left(\sqrt{h}\sin(\theta)\xi_2\right)d\theta\\
&=2\int_{0}^{\frac{\pi}{2}} \cos\left(\sqrt{h}\sin(\theta)\xi_2\right)d\theta\\
&=\pi\mathtt{J}_0(\sqrt{h}\xi_2),
\end{align*}
where $\mathtt{J}_0$ denotes the first Bessel function of order zero. Therefore, equation (\ref{aa}) reads
\begin{align*}
&\quad\mathcal{F}\Gamma(t-s,x_1,\cdot,x_2-\cdot)(\xi)\\
&\quad=\frac{e^{-ix_2\xi_2}}{2}\int_{x_1}^{x_1+2(t-s)} e^{-iy_1\xi_1}\mathtt{J}_0(\sqrt{h}\xi_2)dy_1\\
&\quad=\frac{e^{-ix_2\xi_2}}{2}\int_{x_1}^{x_1+2(t-s)} e^{-iy_1\xi_1}\mathtt{J}_0\left(\xi_2\sqrt{(2(t-s)+x_1-y_1)(y_1^3-x_1^3)/3}\right)dy_1.
\end{align*}
We now set $\tau:=t-s>0$ and for $\lambda\in ]0,1[$ we use the change of variable $y_1=x_1+2\lambda\tau$; this gives
\begin{align}\label{as}
&\mathcal{F}\Gamma(\tau,x_1,\cdot,x_2-\cdot)(\xi)\nonumber\\
&\quad=\tau e^{-ix_2\xi_2}\int_{0}^{1} e^{-i(x_1+2\lambda\tau)\xi_1}\mathtt{J}_0\left(\xi_2\sqrt{4\tau^2\lambda(1-\lambda)(4\tau^2\lambda^2+6\tau x_1\lambda+3x_1^2)/3}\right)d\lambda\nonumber\\
&\quad=\tau e^{-ix_2\xi_2-ix_1\xi_1}\int_{0}^{1} e^{-i2\lambda\tau\xi_1}\mathtt{J}_0\left(\xi_2\sqrt{4\tau^2\lambda(1-\lambda)(4\tau^2\lambda^2+6\tau x_1\lambda+3x_1^2)/3}\right)d\lambda\nonumber\\
&\quad=\tau e^{-ix_2\xi_2-ix_1\xi_1}\int_{0}^{1} e^{-i2\lambda\tau\xi_1}\mathtt{J}_0\left(\xi_2\tilde{h}(\tau,x_1,\lambda)\right)d\lambda,
\end{align}
where we introduced the shorthand notation
\begin{align}\label{h tilde}
\tilde{h}(\tau,x_1,\lambda):=&\sqrt{4\tau^2\lambda(1-\lambda)(4\tau^2\lambda^2+6\tau x_1\lambda+3x_1^2)/3}\nonumber\\
=& 2\tau\sqrt{\lambda(1-\lambda)}\sqrt{(4\tau^2\lambda^2+6\tau x_1\lambda+3x_1^2)/3}.
\end{align}
Taking the modulus of the first and last members in (\ref{as}) we see that
\begin{align}\label{aq}
|\mathcal{F}\Gamma(\tau,x_1,\cdot,x_2-\cdot)(\xi)|\leq\tau\left|\int_{0}^{1} e^{-i2\lambda\tau\xi_1}\mathtt{J}_0\left(\xi_2\tilde{h}(\tau,x_1,\lambda)\right)d\lambda\right|.
\end{align}

To establish the bound \eqref{condition for integral}, we can focus on the behaviour of $\xi\mapsto|\mathcal{F}\Gamma(\tau,x_1,\cdot,x_2-\cdot)(\xi)|$ for large values of $|\xi|$ only; in fact, such function is smooth and bounded on any compact set containing the origin. According to formula (1) page 206 in \cite{Watson}, the Bessel function $\mathtt{J}_0$ can be represented for $z>0$ as
\begin{align}\label{Bessel repr}
\mathtt{J}_0(z)=\sqrt{\frac{2}{\pi z}}\left[\cos\left(z-\frac{\pi}{4}\right)\mathtt{P}_+(z)-\sin\left(z-\frac{\pi}{4}\right)\mathtt{P}_-(z)\right]
\end{align}
with
\begin{align*}
\mathtt{P}_{+}(z)=\frac{1}{2\sqrt{\pi}}\int_0^{+\infty}e^{-u}\frac{1}{\sqrt{u}}\left\{\left(1+\frac{iu}{2z}\right)^{-\frac{1}{2}}+\left(1-\frac{iu}{2z}\right)^{-\frac{1}{2}}\right\}du
\end{align*}
and
\begin{align*}
\mathtt{P}_{-}(z)=\frac{1}{2i\sqrt{\pi}}\int_0^{+\infty}e^{-u}\frac{1}{\sqrt{u}}\left\{\left(1+\frac{iu}{2z}\right)^{-\frac{1}{2}}-\left(1-\frac{iu}{2z}\right)^{-\frac{1}{2}}\right\}du.
\end{align*}
Since $\left|1+\frac{iu}{2z}\right|=\sqrt{1+\frac{u^2}{4z^2}}\geq 1$, we see that $|\mathtt{P}_{\pm}(z)|\leq 1$; this together with \eqref{Bessel repr} implies $|\mathtt{J}_0(z)|\leq\frac{2\sqrt{2}}{\sqrt{\pi}}\frac{1}{\sqrt{z}}$ and  
\begin{align}\label{first bound}
\left|\int_{0}^{1} e^{-i2\lambda\tau\xi_1}\mathtt{J}_0\left(\xi_2\tilde{h}(\tau,x_1,\lambda)\right)d\lambda\right|\leq& C_1\int_0^1\frac{1}{|\xi_2|^{\frac{1}{2}}|\tilde{h}(\tau,x_1,\lambda)|^{\frac{1}{2}}}d\lambda\nonumber\\
=&\frac{C_1}{|\xi_2|^{\frac{1}{2}}}\int_0^1\frac{1}{|\tilde{h}(\tau,x_1,\lambda)|^{\frac{1}{2}}}d\lambda.
\end{align}
We now evaluate the last integral above; recalling the definition of $\tilde{h}$ in \eqref{h tilde} we can write
\begin{align*}
\int_0^1\frac{1}{|\tilde{h}(\tau,x_1,\lambda)|^{\frac{1}{2}}}d\lambda=\int_0^1\frac{1}{\sqrt{2\tau}}\frac{1}{(\lambda(1-\lambda))^{\frac{1}{4}}}\frac{1}{L^{\frac{1}{4}}}d\lambda,
\end{align*}
where $L:=(4\tau^2\lambda^2+6\tau x_1\lambda+3x_1^2)/3\geq \tau^2\lambda^2/3$ and hence $L^{\frac{1}{4}}\geq\sqrt{\tau}\sqrt{\lambda}/3^{\frac{1}{4}}$. This gives
\begin{align*}
\int_0^1\frac{1}{|\tilde{h}(\tau,x_1,\lambda)|^{\frac{1}{2}}}d\lambda\leq \frac{C_2}{\tau}\int_0^1\frac{1}{\lambda^{3/4}(1-\lambda)^{1/4}}=\frac{\tilde{C}_2}{\tau}.
\end{align*}
A combination of this estimate with (\ref{first bound}) and (\ref{aq}) yields
\begin{align}\label{beta}
|\mathcal{F}\Gamma(\tau,x_1,\cdot,x_2-\cdot)(\xi)|\leq \frac{C}{|\xi_2|^{\frac{1}{2}}}.
\end{align}
To get a bound involving also $|\xi_1|$, we go back to the integral in (\ref{aq}) and perform an integration by parts (recall from \eqref{h tilde} that $\tilde{h}(\tau,x_1,0)=\tilde{h}(\tau,x_1,1)=0$ and $\mathtt{J}_0(0)=0$); this yields
\begin{align*}
&\int_{0}^{1} e^{-i2\lambda\tau\xi_1}\mathtt{J}_0\left(\xi_2\tilde{h}(\tau,x_1,\lambda)\right)d\lambda\\
&\quad=\frac{1}{i2\tau\xi_1}-\frac{e^{-i2\tau\xi_1}}{i2\tau\xi_1}+\frac{1}{i2\tau\xi_1}\int_0^1e^{-i2\lambda\tau\xi_1}\mathtt{J}_0'\left(\xi_2\tilde{h}(\tau,x_1,\lambda)\right)\xi_2\partial_{\lambda}\tilde{h}(\tau,x_1,\lambda)d\lambda. 
\end{align*}
Recalling that for all $z \in\mathbb{C}$ we have $\mathtt{J}_0'(z)=-\mathtt{J}_1(z)$, the Bessel function of order one, and that $|\mathtt{J}_1(x)|\leq M$, for all $x\in\mathbb{R}$ and a suitable positive constant $M$, we can write 
\begin{align}\label{second bound}
&\left|\int_{0}^{1} e^{-i2\lambda\tau\xi_1}\mathtt{J}_0\left(\xi_2\tilde{h}(\tau,x_1,\lambda)\right)d\lambda\right|\nonumber\\
&\quad\leq \frac{C_4}{\tau|\xi_1|}+\frac{C_5|\xi_2|}{\tau|\xi_1|}\int_0^1\left|\mathtt{J}_0'\left(\xi_2\tilde{h}(\tau,x_1,\lambda)\right)\right||\partial_{\lambda}\tilde{h}(\tau,x_1,\lambda)|d\lambda\nonumber\\
&\quad\leq\frac{C_4}{\tau|\xi_1|}+\frac{MC_5|\xi_2|}{\tau|\xi_1|}\int_0^1|\partial_{\lambda}\tilde{h}(\tau,x_1,\lambda)|d\lambda. 
\end{align}
Now,
\begin{align*}
\int_0^1|\partial_{\lambda}\tilde{h}(\tau,x_1,\lambda)|d\lambda=&\tau\int_0^1\frac{|1-2\lambda|}{\sqrt{\lambda(1-\lambda)}}L^{1/2}d\lambda\\
&+\tau\int_0^1\frac{\sqrt{\lambda(1-\lambda)}}{|L^{-1/2}|}|\partial_{\lambda}L|d\lambda\\
\leq&c_1\sqrt{\tau^2+x_1^2}+c_2(\tau+|x_1|)\\
\leq& c_3(\tau+|x_1|);
\end{align*}
here, we utilized in the first integral the bound $L\leq C(\tau^2+x_1^2)$ while in the second $L^{1/2}\geq \frac{1}{\sqrt{3}}\tau\lambda$. Therefore, using this estimate with \eqref{second bound} in \eqref{aq} we get
\begin{align}\label{alpha}
|\mathcal{F}\Gamma(\tau,x_1,\cdot,x_2-\cdot)(\xi)|&\leq\frac{C_4}{|\xi_1|}+\frac{MC_5|\xi_2|}{|\xi_1|}c_3(\tau+|x_1|)\nonumber\\
&=\frac{C_4}{|\xi_1|}+K(\tau,x_1)\frac{|\xi_2|}{|\xi_1|}
\end{align}
Combining (\ref{beta}) and (\ref{alpha}) we can now complete the estimate of $|\mathcal{F}\Gamma(\tau,x_1,\cdot,x_2-\cdot)(\xi)|$ for large value of $|\xi|$. Let $\theta\in]0,1[$ (to be fixed later):
\begin{itemize}
\item if $\xi$ is such that $|\xi_2|\leq M|\xi_1|^{\theta}$, for some positive constant $K$, then from inequality (\ref{alpha}) we get
\begin{align}\label{1}
|\mathcal{F}\Gamma(\tau,x_1,\cdot,x_2-\cdot)(\xi)|&\leq \frac{C_4}{|\xi_1|}+K(\tau,x_1)M\frac{|\xi_1|^{\theta}}{|\xi_1|}\nonumber\\
&=\frac{C_4}{|\xi_1|}+K(\tau,x_1)M\frac{1}{|\xi_1|^{1-\theta}}\nonumber\\
&\leq \tilde{K}(\tau,x_1)\frac{1}{|\xi_1|^{1-\theta}}\nonumber\\
&\leq \tilde{K}(\tau,x_1)\frac{1}{|\xi|^{1-\theta}};
\end{align}
here, the last inequality is due to $|\xi|\lesssim |\xi_1|+|\xi_2|\lesssim |\xi_1|+|\xi_1|^{\theta}\lesssim |\xi_1|$ and hence $\frac{1}{|\xi_1|}\lesssim\frac{1}{|\xi|}$;
\item otherwise, if $\xi$ is such that $|\xi_2|\geq M|\xi_1|^{\theta}$, then 
$|\xi|\lesssim |\xi_1|+|\xi_2|\lesssim |\xi_2|^{1/\theta}+|\xi_2|\lesssim |\xi_2|^{1/\theta}$ or equivalently $|\xi|^{\theta}\lesssim |\xi_2|$. This condition, combined with (\ref{beta}), yields
\begin{align}\label{2}
|\mathcal{F}\Gamma(\tau,x_1,\cdot,x_2-\cdot)(\xi)|&\leq\frac{C}{|\xi|^{\theta/2}}.
\end{align}
\end{itemize}
To match the exponents in (\ref{1}) and (\ref{2}) we need to impose $1-\theta=\theta/2$, that means $\theta=2/3$. Hence, for large values of $|\xi|$ we get the estimate
\begin{align}\label{final}
|\mathcal{F}\Gamma(\tau,x_1,\cdot,x_2-\cdot)(\xi)|^2&\leq \frac{\kappa(\tau,x_1)^2}{|\xi|^{2/3}},
\end{align}
with $\kappa(\tau,x_1)$ being with linear growth in $\tau$ and $|x_1|$. For a global (in $\xi$) estimate we can simply set
\begin{align}\label{final}
|\mathcal{F}\Gamma(\tau,x_1,\cdot,x_2-\cdot)(\xi)|^2&\leq \frac{\tilde{\kappa}(\tau,x_1)}{1+|\xi|^{2/3}},\quad\xi\in\mathbb{R}^2,
\end{align}
and conclude that
\begin{align*}
&\int_0^t\int_{\mathbb{R}^2}|\mathcal{F}\Gamma(t-s,x_1,\cdot,x_2-\cdot)(\xi_1,\xi_2)|^2\mu(d\xi_1d\xi_2)ds\\
&\quad\leq\int_0^t\tilde{\kappa}(t-s,x_1)\int_{\mathbb{R}^2}\frac{1}{1+|\xi|^{2/3}}\mu(d\xi_1d\xi_2)ds,
\end{align*}
which turns out to be finite by virtue of (\ref{SC}) (and the nice behaviour of $\tilde{\kappa}(t-s,x_1)$ with respect to its first argument).

We now prove that the map $(t,x_1,x_2)\mapsto u(t,x_1,x_2)$ is measurable by showing the $\mathbb{L}^2(\Omega)$-continuity of $\{u(t,x)\}_{t\geq 0,x\in\mathbb{R}^2}$. Starting with the time increment we can write for $t\in [0,T]$, $x\in\mathbb{R}^2$ and $h>0$ that
\begin{align*}
&\mathbb{E}\left[|u(t,x_1,x_2)-u(t+h,x_1,x_2)|^2\right]\\
&\quad\leq 2\mathbb{E}\left[\left|\int_0^t\int_{\mathbb{R}^2}\Gamma(t-s,x_1,y_1,x_2-y_2)-\Gamma(t+h-s,x_1,y_1,x_2-y_2)M(ds,dy_1,dy_2)\right|^2\right]\\
&\quad\quad+2\mathbb{E}\left[\left|\int_t^{t+h}\int_{\mathbb{R}^2}\Gamma(t+h-s,x_1,y_1,x_2-y_2)M(ds,dy_1,dy_2)\right|^2\right]\\
&\quad= 2\int_0^t\int_{\mathbb{R}^2}\left|(\mathcal{F}\Gamma(t-s,x_1,\cdot,x_2-\cdot))(\xi)-(\mathcal{F}\Gamma(t+h-s,x_1,\cdot,x_2-\cdot)(\xi)\right|^2d\mu(\xi)ds\\
&\quad\quad+2\int_t^{t+h}\int_{\mathbb{R}^2}\left|(\mathcal{F}\Gamma(t+h-s,x_1,\cdot,x_2-\cdot))(\xi)\right|^2d\mu(\xi)ds\\
&\quad= 2\int_0^t\int_{\mathbb{R}^2}\left|(\mathcal{F}\Gamma(t-s,x_1,\cdot,\cdot))(\xi)-(\mathcal{F}\Gamma(t+h-s,x_1,\cdot,\cdot)(\xi)\right|^2d\mu(\xi)ds\\
&\quad\quad+2\int_t^{t+h}\int_{\mathbb{R}^2}\left|(\mathcal{F}\Gamma(t+h-s,x_1,\cdot,\cdot))(\xi)\right|^2d\mu(\xi)ds.
\end{align*}
The integrand in the first integral is bounded by a constant times $\frac{1}{1+|\xi|^{2/3}}$, which is by assumption integrable with respect to $\mu$; therefore, by dominated convergence the first integral tends to zero as $h\to 0$. The second integral also converges to zero by virtue of (\ref{condition for integral}). \\
The increment in the variable $x_2$ is treated similarly; in fact, for $t\in [0,T]$, $x_1,x_2,z_2\in\mathbb{R}$ we have
\begin{align*}
&\mathbb{E}\left[|u(t,x_1,x_2)-u(t,x_1,z_2)|^2\right]\\
&\quad= \int_0^t\int_{\mathbb{R}^2}\left|(\mathcal{F}\Gamma(t-s,x_1,\cdot,x_2-\cdot))(\xi)-(\mathcal{F}\Gamma(t-s,x_1,\cdot,z_2-\cdot)(\xi)\right|^2d\mu(\xi)ds\\
&\quad=\int_0^t\int_{\mathbb{R}^2}|e^{ix_2\xi_2}-e^{iz_2\xi_2}|^2\left|(\mathcal{F}\Gamma(t-s,x_1,\cdot,\cdot))(\xi)\right|^2d\mu(\xi)ds,
\end{align*}
and by dominated convergence we conclude that the last integral tends to zero as $|x_2-z_2|\to 0$. Lastly, we consider the increment in the variable $x_1$: for  $t\in [0,T]$, $x_1,z_1,x_2\in\mathbb{R}$ we get
\begin{align*}
&\mathbb{E}\left[|u(t,x_1,x_2)-u(t,z_1,x_2)|^2\right]\\
&\quad= \int_0^t\int_{\mathbb{R}^2}\left|(\mathcal{F}\Gamma(t-s,x_1,\cdot,x_2-\cdot))(\xi)-(\mathcal{F}\Gamma(t-s,z_1,\cdot,x_2-\cdot)(\xi)\right|^2d\mu(\xi)ds\\
&\quad= \int_0^t\int_{\mathbb{R}^2}\left|(\mathcal{F}\Gamma(t-s,x_1,\cdot,\cdot))(\xi)-(\mathcal{F}\Gamma(t-s,z_1,\cdot,\cdot)(\xi)\right|^2d\mu(\xi)ds.
\end{align*}
The integrand above can be upper bounded by a constant, which depends on $x_1$ and $z_1$, times $\frac{1}{1+|\xi|^{2/3}}$; this fact together with dominated convergence implies that the last integral tends to zero as $|x_1-z_1|\to 0$.

\end{document}